\documentclass[12pt,a4paper]{amsart}
\usepackage[utf8]{inputenc}
\usepackage{amsfonts}
\usepackage{amssymb}
\usepackage{amsmath}
\usepackage{amsthm}
\usepackage{cite}
\usepackage{enumitem}
\usepackage{tikz}
\usepackage{subfigure}

\theoremstyle{plain}
\newtheorem{thm}{Theorem}  
\newtheorem{defn}{Definition}

\newtheorem{rem}{Remark} 
\newtheorem{conj}{Conjecture}

\providecommand{\sm}{\setminus}
\providecommand{\N}{\mathbb{N}}
\providecommand{\R}{\mathbb{R}}
\providecommand{\Z}{\mathbb{Z}}

\providecommand{\ov}{\overline}

\DeclareMathOperator{\dist}{dist}

\renewcommand{\qed}{\hfill $\Box$}

\setitemize{itemsep=+2pt}
\setenumerate{itemsep=+2pt}
\def\XXint#1#2#3{{\setbox0=\hbox{$#1{#2#3}{\int}$} 
     \vcenter{\hbox{$#2#3$}}\kern-.5\wd0}}

\begin{document}

\allowdisplaybreaks

\title{A uniqueness result for the Sine-Gordon breather}

\author{Rainer Mandel}
\address{R. Mandel\hfill\break
Karlsruhe Institute of Technology \hfill\break
Institute for Analysis \hfill\break
Englerstra{\ss}e 2 \hfill\break
D-76131 Karlsruhe, Germany}
\email{rainer.mandel@kit.edu}
\date{\today}


\begin{abstract}
   In this note we  prove that the sine-Gordon breather is the only quasimonochromatic breather in the context
   of nonlinear wave equations in $\R^N$.
\end{abstract}

\maketitle

\section{Introduction}

  Breathers are time-periodic and spatially localized patterns that describe the propagation of waves. The
  most impressive solution of this kind is the so-called sine-Gordon breather for the 1D
  sine-Gordon equation 
  $$
    \partial_{tt}u - \partial_{xx} u + \sin(u) = 0 \qquad\text{in }\R\times\R.
  $$ 
  It is given by the explicit formula  
  \begin{align} \label{eq:sineGordon}
    u^*(x,t) = 4\arctan\left(\frac{m\sin(\omega t)}{\omega\cosh(mx)} \right) 
    \qquad\text{for }(x,t)\in\R\times\R,
  \end{align} 
  where the parameters $m,\omega>0$ satisfy $m^2+\omega^2 = 1$. It is natural to ask if other real-valued
  breather solutions exist. We shall address this question in the broader context of more general nonlinear
  wave equations of the form
  \begin{equation}\label{eq:NLWave}
    \partial_{tt}u- \Delta u = g(u) \quad\text{in }\R^N\times\R,
  \end{equation} 
  where the space dimension $N\in\N$ and the nonlinearity $g:\R\to\R$ are arbitrary. 
  
  \medskip
  
  The existence of  radially symmetric breather solutions for the cubic Klein-Gordon equation $g(z)=-m^2z +
  z^3,\, m>0$ in three spatial dimensions was established in~\cite{Scheider_Breather}. These real-valued 
  solutions are only weakly localized in the sense that they satisfy $u(\cdot,t)\in L^q(\R^N)$ for some $q\in
  (2,\infty)$ but  $u(\cdot,t)\notin L^2(\R^N)$. In~\cite{ManSch_Variational} infinitely
  many weakly localized breathers were found for nonlinearities  $Q(x)|u|^{p-2}u$ where $Q$ lies in a
  suitable Lebesgue space and $p>2$ is chosen suitably depending on $Q$ as well as the space dimension $N\geq 2$. Up
  to now, nothing is known about the existence of strongly localized breathers of~\eqref{eq:NLWave}
  satisfying $u(\cdot,t)\in L^2(\R^N)$ for almost all $t\in\R$ and $N\geq 2$, see
  however~\cite{PlumReichel_Breather} for a an existence result for semilinear curl-curl equations for $N=3$.
  In the  case $N=1$ strongly localized breather solutions different from the sine-Gordon breather
  have been found for nonlinear wave equations of the form 
  $$
    s(x)\partial_{tt}u - u_{xx} + q(x)u = f(x,u) \qquad (x\in\R)
  $$ 
  where the coefficient functions $s,q$ are discontinuous and periodic, see~\cite[Theorem~1.3]{Hirsch} and
  ~\cite[Theorem~1.1]{Schneider}. Given the discontinuity of $s,q$ it must be expected that these breathers
  are not twice continuously differentiable. To sum up, the existence of smooth and strongly localized breather
  solutions of~\eqref{eq:NLWave} different from the sine-Gordon breather is not known. Still for $N=1$
  there are   nonexistence results by Denzler~\cite{Denzler_Nonpersistence} and Kowalczyk,
  Martel, Mu\~{n}oz \cite{KowMarMun_Nonexistence} dealing with small perturbations of the sine-Gordon
  equation respectively small odd breathers (not covering the even sine-Gordon breather). We are not aware of
  any other mathematically rigorous existence or nonexistence results for~\eqref{eq:NLWave}.

  \medskip
  
  One of the main obstructions for the construction of localized breathers is polychromaticity. Indeed,
  plugging in an ansatz of the form $u(x,t)=\sum_{k\in\Z} u_k(x)e^{ikt}$ with $u_k=\ov{u_{-k}}$ one ends up
  with infinitely many equations of nonlinear Helmholtz type that typically do not possess strongly localized
  solutions, see for instance~\cite[Theorem~1a]{Kato_Growth}. For this reason the solutions obtained in  
  \cite{Scheider_Breather,ManSch_Variational} are only weakly localized. On the other hand, a
  purely monochromatic ansatz like $u(x,t)=\sin(\omega t)p(x)$ cannot be successful either provided that $g$
  is not a linear function. In view of the formula~\eqref{eq:sineGordon} for the sine-Gordon breather we
  investigate whether quasimonochromatic breathers exist.
    
  \begin{defn} \label{defn}
    We call the function $u:\R^N\times\R\to\R$  a quasimonochromatic breather if  
    $$
      u(x,t)=F(\sin(\omega t)p(x)) \qquad (x\in\R^N,t\in\R)
    $$ 
    for some $\omega\in\R\sm\{0\}$ and nontrivial functions $F\in C^2(\R),p\in C^2(\R^N)$ such that 
    $F(0)=0$ and $p(x)\to 0$ as $|x|\to\infty$. 
  \end{defn}
  
  We show that in one spatial dimension the sine-Gordon breather is, up to translation and dilation,
  the only one for~\eqref{eq:NLWave} and that no such breathers exist in higher dimensions as long as $g$
  does not act like a linear function. In fact, to rule out $L^\infty$-small solutions of linear wave
  equations, we assume that $g:\R\to\R$ is not a linear function near zero, i.e., that there is a nontrivial
  interval $I\subset\R$ containing $0$ with the property that there is no $\beta\in\R$ such that $g(z)=\beta z$
  for all $z\in I$. 
  
  \begin{thm} \label{thm:Nonexistence}
    Assume $N\in\N$ and that $g:\R\to\R$ is not a linear function near zero.   
    \begin{itemize}
      \item[(i)] In the case $N\geq 2$ there is no quasimonochromatic breather solution of~\eqref{eq:NLWave}.
      \item[(ii)] In the case $N=1$ each quasimonochromatic breather solution of~\eqref{eq:NLWave} is of the
      form $u(x,t)= \kappa u^*(x-x_0,t)$ for $x_0\in\R$, $m,\omega,\kappa\in\R\sm\{0\}$ and $u^*$ as
      in~\eqref{eq:sineGordon}. The nonlinearity then satisfies $g(z) =  -(m^2+\omega^2)\kappa
      \sin(\kappa^{-1} z)$ whenever $|z|< 2\pi|\kappa|$.
    \end{itemize}
  \end{thm}
  
  We stress that our result holds regardless of any smoothness assumption on $g$ nor any kind of growth
  condition at 0 or infinity. Moreover, our considerations are not limited to small perturbations of $u^*$ or
  small breathers in whatever sense. Following the proof of Theorem~\ref{thm:Nonexistence} one also finds that
  quasimonochromatic breathers of wave equations on any open set  $\Omega\subsetneq\R^N$ with homogeneous
  Dirichlet conditions 
  \begin{equation} \label{eq:waveOmega}
    \partial_{tt}u- \Delta u = g(u) \quad\text{in }\Omega\times\R,\qquad
    u=0 \text{ on }\partial\Omega\times\R 
  \end{equation}
  with profile functions $p\in C^2(\ov\Omega)$ do not exist either (even if $N=1$) provided that $g$ is not
  a linear function near zero. We will comment on this fact at the end of this paper. As a consequence,
  we find that Rabinowitz' $C^2([0,1]\times \R)$-solutions of the 1D wave equation from~\cite[Theorem~1.6]{Rab_vibes} are
  not of quasimonochromatic type. This might be true as well for the solutions
  from\cite{Brezis_vibes,Coron_Periodic}, but here our argument does not apply in a
  direct way since the solutions are not known to be twice continuously differentiable up to the
  boundary.

  \medskip
  
  For completeness we briefly comment on the linear case  $g(z)=\beta z$, $\beta\in\R$.
  Then the profile function $p$ of any given quasimonochromatic breather of~\eqref{eq:NLWave} satisfies the
  linear elliptic PDE $ -\Delta p - (\omega^2+\beta)p  = 0$ in $\R^N$. 
  For $\beta<-\omega^2$ there are positive, radially symmetric and exponentially decaying solutions $p$, see~\cite[Theorem~2]{GNN_Symmetry}. 
  In the case $\beta>-\omega^2,N\geq 2$ one can find radial as well as non-radial solutions of the
  associated Helmholtz equation all of which have infinitely many nodal domains and satisfy $|p(x)|+|\nabla p(x)|\gtrsim
  |x|^{\frac{1-N}{2}}$ in a suitable integrated sense,
  see~\cite[Theorem~1]{Swanson_semilinear} respectively~\cite[Theorem~1a]{Kato_Growth}. For
  $\beta>-\omega^2,N=1$ all solutions are linear combinations of $\sin$ and $\cos$ so that breather solutions
  do not exist. So we see that the picture is already quite complete in the case of linear wave equations.

\section{Proof of Theorem~\ref{thm:Nonexistence}}

   In the following let $u(x,t)=F(\sin(\omega t)p(x))$ be a solution of~\eqref{eq:NLWave} as in~\eqref{defn}
   with $g$ as in the Theorem.  
    Plugging in this ansatz we get for all $x\in\R^N$ such that $p(x)\neq 0$,
  \begin{align*}
    \partial_{tt}u(x,t)
    &= - \omega^2\sin(\omega t)p(x)F'(\sin(\omega t)p(x)) + \omega^2\cos(\omega t)^2p(x)^2 F''(\sin(\omega
    t)p(x)) \\
    &= - \omega^2 zF'(z) + \omega^2(p(x)^2-z^2) F''(z), \\
    \Delta u(x,t) 
    &= \sin(\omega t)\Delta p(x) F'(\sin(\omega t)p(x)) +  \sin(\omega t)^2|\nabla p(x)|^2F''(\sin(\omega
    t)p(x)) \\
    &=   \frac{\Delta p(x)}{p(x)} zF'(z) +  \frac{|\nabla p(x)|^2}{p(x)^2} z^2 F''(z),
  \end{align*}
  where $z=\sin(\omega t)p(x)\in [-\|p\|_\infty,+\|p\|_\infty]$.  
  This and~\eqref{eq:NLWave} imply for $x\in\R^N,z\in\R$ such that $p(x)\neq 0, z\in [-\|p\|_\infty,+\|p\|_\infty]$
  \begin{align} \label{eq:AlmostReducedEquations}
  \begin{aligned}
     g(F(z))+\omega^2 z F'(z) + \omega^2z^2F''(z) 
     =    p(x)^2 \omega^2 F''(z)-  \frac{\Delta p(x)}{p(x)} z F'(z) -  \frac{|\nabla
    p(x)|^2}{p(x)^2} z^2 F''(z).
    \end{aligned} 
  \end{align}
  If $F$ was linear on $[-\|p\|_\infty,+\|p\|_\infty]$, then $g$ would have to be linear 
  on the nontrivial interval $I:=\{F(z):|z|\leq \|p\|_\infty\}$ as well. Since the latter is not the case by
  assumption, we know that $z\mapsto z^2 F''(z)$ does not vanish identically on that interval.
  Multiplying~\eqref{eq:AlmostReducedEquations} with $p(x)$ and choosing $z$ according to 
  $z^2 F''(z)\neq 0$ we
  find that $p$ does not change sign. Indeed, if $p(x^*)\neq 0$ and $R>0$ is the smallest radius such
  that $p$ has a fixed sign in the open ball $B_R(x^*)$, then  Hopf's Lemma
  \cite[Lemma~3.4]{GiTr} implies $|\nabla p|>0$ on $\partial B_R(x^*)$. But
  then~\eqref{eq:AlmostReducedEquations} implies that $\Delta p$ is unbounded on $\partial B_R(x^*)$, which
  contradicts $p\in C^2(\R^N)$. Hence, $p$ does not change sign and we will without loss of generality assume
  that $p$ is positive. So~\eqref{eq:AlmostReducedEquations} holds for all $x\in\R^N$ and all $z\in
  [-\|p\|_\infty,\|p\|_\infty]$ and standard elliptic regularity theory gives $p\in C^\infty(\R^N)$.

  \medskip
  
  Differentiating~\eqref{eq:AlmostReducedEquations} with respect to $x_i$ we get  
  \begin{equation} \label{eq:ReductionI}
     \partial_i(p(x)^2) \omega^2F''(z)-  \partial_i\left(\frac{\Delta p(x)}{p(x)}\right) z F'(z) -
    \partial_i \left(\frac{|\nabla p(x)|^2}{p(x)^2}\right) z^2 F''(z) = 0. 
  \end{equation}
  Since $p^2$ is non-constant, we infer that  $F$ satisfies an ODE of the form 
  \begin{equation} \label{eq:ReducedEqI}
    F''(z) = \frac{-\mu_2 z}{\omega^2+\mu_1z^2}  F'(z)
    \qquad (|z|\leq \|p\|_\infty,\; \mu_1\in\R,\mu_2\in\R\sm\{0\}).
  \end{equation}
  Here, $\mu_2\neq 0$ is due to the fact that $F$ is not a linear function.  
  Each nontrivial solution of such an ODE satisfies $F'(z)\neq 0$  for almost all $z\in
  [-\|p\|_\infty,\|p\|_\infty]$. Combining \eqref{eq:ReductionI} and~\eqref{eq:ReducedEqI} we thus infer
  $$
   -\partial_i(p(x)^2)  \frac{\mu_2  \omega^2  z}{\omega^2+\mu_1z^2}    -  \partial_i\left(\frac{\Delta
    p(x)}{p(x)}\right) z   
    + \partial_i \left(\frac{|\nabla p(x)|^2}{p(x)^2}\right) \frac{\mu_2
    z^3}{\omega^2+\mu_1z^2}   = 0. 
  $$
  Since~\eqref{eq:ReducedEqI} holds for all $i\in\{1,\ldots,N\}$ and $z\in [-\|p\|_\infty,\|p\|_\infty]$, we get
  \begin{align*}
    -\mu_1    \partial_i\left(\frac{\Delta
    p(x)}{p(x)}\right) + \mu_2 \partial_i \left(\frac{|\nabla p(x)|^2}{p(x)^2}\right) 
    &= 0, \\ 
    - \mu_2 \partial_i(p(x)^2)    -  \partial_i\left(\frac{\Delta
    p(x)}{p(x)}\right) &= 0.
  \end{align*}
  Since $\mu_2\neq 0$ we can find $\lambda_1,\lambda_2\in\R$ such that  
  $$
    -\mu_1   \frac{\Delta p}{p}  + \mu_2  \frac{|\nabla p|^2}{p^2}
    = - \lambda_2\mu_1+\lambda_1\mu_2,\qquad    
    - \mu_2 p^2  -   \frac{\Delta  p }{p } = - \lambda_2. 
  $$
  This implies
  \begin{align} \label{eq:ReducedEqII}
     |\nabla  p|^2  =  \lambda_1 p^2   - \mu_1   p^4,\qquad 
    -   \Delta p  + \lambda_2 p = \mu_2 p^3. 
  \end{align}

  \medskip
  
  We now use~\eqref{eq:ReducedEqII} and the positivity of $p$ to show that $p$ is radially symmetric about
  its maximum point $x_0\in\R^N$. We concentrate on the case $N\geq 2$ since the claim for $N=1$ follows from the fact that
  $x\mapsto u(x_0+x)$ and $x\mapsto u(x_0-x)$ solve the same initial value
  problem. Since $p$ vanishes at infinity, we must have $\lambda_1\geq 0$ and, since $p$ does not change sign,
  $\lambda_2\geq 0$, see~\cite[Theorem~1]{Swanson_semilinear}. Moreover, $p$ attains its maximum at some point
  $x_0\in\R^N$ with $p(x_0)>0, |\nabla p(x_0)|=0,\Delta p(x_0)\leq 0$. This and~\eqref{eq:ReducedEqII} implies $\lambda_1,\mu_1>0$ as well
  as $\mu_2\geq 0$. So we know that \eqref{eq:ReducedEqII} holds for 
  $$
    \lambda_1,\mu_1>0,\qquad \lambda_2,\mu_2\geq 0.
  $$
  In the case $\lambda_2>0$ Theorem~2 from~\cite{GNN_Symmetry} implies the radial
  symmetry about $x_0$, so we are left with the case $\lambda_2=0$. 
  
  \medskip
  
  So let use assume $\lambda_2=0$. Liouville's Theorem implies that $\mu_2=0$ is impossible, so we have
  $\mu_2>0$ in this case. Define $\alpha:=1-\frac{\mu_2}{\mu_1}\in (-\infty,1)$.
  In the case $\alpha\in (0,1)$ the function $\psi(x):= p(x)^\alpha$ satisfies 
  \begin{align*}
    -\Delta \psi 
    = - \alpha (\Delta p) p^{\alpha-1} - \alpha(\alpha-1) |\nabla p|^2 p^{\alpha-2} 
    \stackrel{\eqref{eq:ReducedEqII}}=  \alpha(1-\alpha)\lambda_1  \psi.  
  \end{align*}
  In view of $\alpha(1-\alpha)\lambda_1>0$ Theorem~1 from~\cite{Swanson_semilinear} implies that $\psi$ has
  infinitely many nodal domains, which contradicts the positivity of $\psi$. So this case cannot occur. In the case $\alpha\in (-\infty,0)$
  radial symmetry about $x_0$ follows once more from~\cite[Theorem~2]{GNN_Symmetry}, so it remains to discuss
  the case $\alpha=0$, i.e., $\mu_1=\mu_2$. Then $\psi(x):=\log(p(x))$ satisfies
  \begin{align*}
    -\Delta \psi 
    = -  (\Delta p) p^{-1} +  |\nabla p|^2 p^{-2} 
     \stackrel{\eqref{eq:ReducedEqII}}=  \lambda_1  \psi  
  \end{align*}
  and we find as above that  $\psi$ has to change sign, which is a contradiction. So we have shown that $p$ is
  radially symmetric about $x_0$ also in the case $\lambda_2=0$.
  
  \medskip
  
  So we have 
  $$
    p(x)=p_0(|x-x_0|)\qquad\text{where } p_0'(r)^2 = \lambda_1 p_0(r)^2 - \mu_1 p_0(r)^4,\quad  p_0'(0)=0.
  $$ 
  Solving this ODE gives 
  $$
    p_0(r) = \frac{A}{\cosh(mr)}\qquad\text{where } \lambda_1= m^2, \mu_1 = m^2A^{-2}
  $$ 
  for some $A>0,m\neq 0$. So $- \Delta p  + \lambda_2 p  = \mu_2 p^3$ can only hold for $N=1$ as well as 
  $\lambda_2= m^2$, $\mu_2 = 2m^2A^{-2}$. 
  Plugging these values into~\eqref{eq:ReducedEqI} and solving the ODE we get from $F(0)=0, F\not\equiv 0$ 
  $$
    F(z) = 4\kappa \arctan\left( \frac{mz}{A\omega}\right) \qquad \text{for some }\kappa\in\R\sm\{0\}.
  $$
  This implies that the breather solution is given by
  $$
    u(x,t) 
    = F(\sin(\omega t)p(x))
    = F(\sin(\omega t)p_0(|x-x_0|))
    = \kappa u^*(x-x_0,t)
  $$
  for $u^*$ as in~\eqref{eq:sineGordon}. So have proved the nonexistence of such breathers for $N\geq 2$
  from claim (i) and the uniqueness statement from claim (ii). 
  
  \medskip
  
  To see that this solution formula determines
  the nonlinearity $g$, we combine~\eqref{eq:ReducedEqI} and~\eqref{eq:ReducedEqII} to get 
  $$  
    p(x)^2 \omega^2 F''(z)-  \frac{\Delta p(x)}{p(x)} z F'(z) -  \frac{|\nabla
    p(x)|^2}{p(x)^2} z^2 F''(z)
    = \frac{m^2(m^2z^2-A^2\omega^2)}{m^2 z^2+A^2\omega^2} F'(z)z. 
  $$  
  So~\eqref{eq:AlmostReducedEquations} implies 
  \begin{align*}
     g(F(z))
     &= -\omega^2 z F'(z) - \omega^2z^2F''(z) + \frac{m^2(m^2 z^2-A^2\omega^2)}{m^2 z^2+A^2\omega^2}
     F'(z)z  \\
     &=  \frac{(m^2+\omega^2)(m^2z^2-A^2\omega^2)}{m^2z^2+A^2\omega^2} zF'(z) \\
     &=  \frac{4Am\kappa\omega(m^2+\omega^2)(m^2z^2-A^2\omega^2)z}{(m^2z^2+A^2\omega^2)^2}.      
  \end{align*}
  Plugging in $z=\frac{A\omega}{m}\tan(\frac{y}{4\kappa})$ for $|y|<2\pi|\kappa|$ we get $F(z)=y$ and hence
  \begin{align*}
    g(y)
    &= 
    \frac{4A^2\omega^2\kappa(m^2+\omega^2)(
    A^2\omega^2\tan(\frac{y}{4\kappa})^2-A^2\omega^2)\tan(\frac{y}{4\kappa})}{
    (A^2\omega^2\tan(\frac{y}{4\kappa})^2+A^2\omega^2)^2}
    \\
    &= 
    \frac{4\kappa(m^2+\omega^2)(\tan(\frac{y}{4\kappa})^2-1)\tan(\frac{y}{4\kappa})}{(\tan(\frac{y}{4\kappa})^2+1)^2}
    \\
    &=  
    4\kappa(m^2+\omega^2)\left(\sin(\frac{y}{4\kappa})^2-\cos(\frac{y}{4\kappa})^2\right)\sin(\frac{y}{4\kappa})\cos(\frac{y}{4\kappa})
    \\
    &=   -2\kappa  (m^2+\omega^2) \cos(\frac{y}{2\kappa})\sin(\frac{y}{2\kappa})    \\
    &=   -\kappa (m^2+\omega^2) \sin(\frac{y}{\kappa}).  
   \end{align*}
  \qed
  
  \begin{rem} ~
    \begin{itemize}
      \item[(i)] 
    We explain why nonlinear quasimonochromatic breathers of ~\eqref{eq:waveOmega} with profile
    functions $p\in C^2(\ov\Omega)$ do not exist on open sets $\Omega\subsetneq\R^N$. 
    The arguments presented above reveal that any such breather is given by functions $F,p$ as in
    Definition~\ref{defn} such that for all $x\in\Omega, p(x)\neq 0, |z|\leq \|p\|_\infty$ we have as in~\eqref{eq:AlmostReducedEquations}
  \begin{align*}  
     g(F(z))+\omega^2 z F'(z) + \omega^2z^2F''(z) 
     =    p(x)^2 \omega^2 F''(z)-  \frac{\Delta p(x)}{p(x)} z F'(z) -  \frac{|\nabla
    p(x)|^2}{p(x)^2} z^2 F''(z).
  \end{align*}
    Now fix $z\in (-\|p\|_\infty,\|p\|_\infty)$ such that $z^2F''(z)\neq 0$ and 
    choose $x^*\in\Omega$ such that
    $p(x^*)\neq 0$. Let $R>0$  be largest possible such that $|p|$ is positive in the open ball
    $B_R(x^*)\subset \Omega$. By the homogeneous Dirichlet boundary condition, 
    we know  $R\leq \dist(x^*,\partial\Omega)<\infty$ and that $p$ vanishes on $\partial
    B_R(x^*)$. So the same argument as in the above proof (Hopf's Lemma) shows that $|\Delta p|$ is unbounded
    on $B_R(x^*)$, a contradiction. As a consequence, such a profile function cannot exist and we obtain the
    nonexistence of quasimonochromatic  breathers for~\eqref{eq:waveOmega}.
    \item[(ii)] In our proof we did not use the assumption $p(x)\to 0$ as $|x|\to\infty$
    when we proved that $|p|$ is positive. As a consequence, each profile function $p$ of a solution
    $u(x,t)=F(\sin(\omega t)p(x))$ of~\eqref{eq:NLWave} has a fixed sign regardless of its behaviour at
    infinity. Similarly,~\eqref{eq:ReducedEqII} holds without this hypothesis. So we conclude that any
    profile function $p\in C^2(\R^N)$ of a quasimonochromatic breather is a positive solution
    of~\eqref{eq:ReducedEqII} provided that the nonlinearity $g$ is not a linear function on the
    interval $\{F(z): |z|\leq \|p\|_\infty\}$. Notice also that the assumption $F(0)=0$ is not used either.
     \item[(iii)] Our notion of a quasimonochromatic breather does not allow for the solutions
     $u(x,t)=u^*(x_1,t)$ ($x\in\R^N$), which are  localized only with respect to one spatial direction.
     Accordingly, our nonexistence result for $N\geq 2$ is false under the 
     weaker
     requirement
     \begin{equation}\label{eq:p_asymptotics}
       \sup_{x'\in\R^{N-1}} |p(x_1,x')|\to 0\qquad\text{as }x_1\to\infty.
     \end{equation}  
     One may conjecture that the solutions $u(x,t)=u^*(x\cdot \theta,t)$  for $\theta\in
     S^{N-1}\subset\R^N$ are the only quasimonochromatic breathers that are localized in some spatial
     direction.
     This open problem bears some similarity to the   Gibbon's Conjecture or 
     de Giorgi Conjecture about the
     classification of monotone solutions of the Allen-Cahn equation $\Delta u + u = u^3$ in $\R^N$ that we
     recast in our setting below.
    \end{itemize}  
  \end{rem}

  \begin{conj}
    Let $N\in\N,N\geq 2$ and let $p\in C^2(\R^N)$ be a solution
    of~\eqref{eq:ReducedEqII} for some $\lambda_1,\lambda_2,\mu_2,\mu_2\in\R$ that satisfies \eqref{eq:p_asymptotics}.
    Then there are $\gamma,m,z\in\R$ such that
    $$
      p(x)=   \frac{\gamma}{\cosh(m(x_1-z))}.
    $$
  \end{conj}
  \medskip
  \begin{conj}
    Let $N\in\N,N\geq 2$ and let $p\in C^2(\R^N)$ be a  solution
    of~\eqref{eq:ReducedEqII} for some $\lambda_1,\lambda_2,\mu_2,\mu_2\in\R$ that satisfies
    $\partial_1 p(x)x_1<0$ for all $x\in\R^N$ such that $x_1\neq 0$.
    Then there are $\gamma,m>0$ such that
    $$
      p(x)=   \frac{\gamma}{\cosh(mx_1)}.
    $$
  \end{conj}


\section*{Acknowledgements}
 
  Funded by the Deutsche Forschungsgemeinschaft (DFG, German Research Foundation)
- Project-ID 258734477 - SFB 1173.
 
\bibliographystyle{plain}
\bibliography{biblio}

\end{document}